\documentclass{amsart}
\usepackage{amssymb}
\usepackage{citesort}

\newtheorem{thm}{Theorem}[section]
\newtheorem{df}[thm]{Definition}
\newtheorem{satz}[thm]{Theorem}
\newtheorem{lem}[thm]{Lemma}
\newtheorem{lemdef}[thm]{Lemma and Definition}
\newtheorem{cor}[thm]{Corollary}

\theoremstyle{definition}

\newtheorem{bez}[thm]{Notation}

\theoremstyle{remark}

\newtheorem{bem}[thm]{Remark}
\newtheorem{rem}[thm]{Remark}

\def\N{{\mathbb{N}}}
\def\R{{\mathbb{R}}}

\def\bew{{\bf Proof: }}

\def\bewvon#1{{\bf Proof of #1: }}
\def\qed{\hspace*{\fill}$\Box$}
\def\rone{{\R^{n+1}}}
\def\reu{{\rone\setminus\{0\}}}
\def\sub{{\rm sub\,}}

\def\dist{{\rm dist\,}}

\def\Dir{{\rm Dir\,}}
\def\oe{{w.~l.~o.~g.\ }}
\def\OE{{W.~l.~o.~g.\ }}
\def\quada{\hspace*{1em}}
\numberwithin{equation}{section}

\begin{document}

\title{Regularity and Splitting of Directed Minimal Cones}

\author{Oliver~C.~Schn\"urer}
\address{Max Planck Institute for Mathematics in the Sciences,
  Inselstr.\ 22-26, 04103 Leipzig, Germany}
\email{Oliver.Schnuerer@mis.mpg.de}


\date{August 2003, reasearch originally from 1998/1999.}


\begin{abstract}
We show that directed minimal cones in $\rone{}$
which have at most one singularity are
-- besides the trivial cases $\emptyset$, $\rone{}$ -- half spaces.
Using blow-up techniques, this result can be used
to get $C^{1,\lambda}$-regularity for the measure-theoretic
boundary of almost minimal Caccioppoli sets which are
representable as subgraphs in $\R^n$, $n\le8$. This provides a different
method to obtain a result due to De Giorgi. We also prove
a splitting theorem for general directed minimal cones. Such a
cone is of the form $\R^k\times C_0$, where $C_0$ is an undirected
minimal cone or a half-line.
\end{abstract}

\maketitle

\section{Introduction}
Let $C\subset\rone{}$ be a minimal cone with its vertex at the origin,
i.~e.\
$$x\in C,\,\tau>0\quad\Longrightarrow\quad \tau x\in C$$
and
$$\int\limits_\Omega|D\chi_C|\le\int\limits_\Omega|D\chi_F|$$
whenever $\Omega\subset\rone{}$ is an open bounded set and
$C\vartriangle F\Subset\Omega$. 
We denote the characteristic function of a set $E$ by $\chi_E$.
A measurable set
$M\subset\rone{}$ is called directed with respect to
$w\in\rone{}\setminus\{0\}$, if the functions
$$\R\ni t\mapsto D_\rho(M,x+tw)\equiv \omega_{n+1}^{-1}\rho^{-(n+1)}
  \int\limits_{B_\rho(x+tw)}\chi_M
  \equiv\int\limits_{B_\rho(x+tw)}\chi_M /
  \int\limits_{B_\rho}1$$
are monotone decreasing for each $\rho>0$ and $x\in\rone{}$.
Such a set is called directed.
Now let us assume additionally that the cone $C$ is directed, nontrivial
($C\not=\emptyset$ and $C\not=\rone{}$) and has at most one singularity.
In this case we prove a regularity theorem which states that
such a cone is a half space. This result was also obtained
by DeGiorgi in \cite{dg}.
\par
Even for non-directed cones, this result is well-known \cite{mm},
provided that the dimension is small, i.~e.\ $n\le6$.
Let $E\subset\rone{}$ be an almost minimal set and 
$x\in\partial M$. We may assume that $x=0$.
Define for $t>0$
$$E_t:=\{x\in\rone{}:tx\in E\}.$$
For an almost minimal set it is well known that there exists a
sequence $t_i\downarrow0$, such that
$$E_{t_i}\to C \mbox{ in } L^1_{loc}(\rone{})$$
and $C$ is a nontrivial minimal cone.
The regularity of $\partial E$ at $0$ is connected with the
regularity of $\partial E$ at $0$ by the equivalence
$$0\in\partial^\star E\quad\Longleftrightarrow\quad
  0\in\partial^\star C$$
for the reduced boundaries.
Especially if $\partial C$ is a hyperplane through the origin
we get $0\in\partial^\star E$. It is well-known \cite{mm} that
in the case of small dimensions ($n\le6$) this blow-up
technique yields equality of the measure-theoretic
boundary $\partial E$ and the reduced boundary $\partial^\star E$
for almost minimal sets $E\subset \rone{}$, because the only
nontrivial minimal cones in $\rone{}$ are half spaces. In
$\R^8$ there are counterexamples.
\par
An application of the regularity theorem is as follows:
The regularity theorem allows us to conclude in the same way
$\partial E=\partial^\star E$ in $\R^8$, when $E$ is a
directed set. This occurs for example when $E$ is representable
as a subgraph.
Let $w\in\R^8\setminus\{0\}$ be the corresponding direction.
Obviously the sets $E_{t_i}$ are again directed with respect to
the same vector $w$. Since the limit in $L^1_{loc}(\rone{})$ of sets which
are directed with respect to $w$ is directed, too, the
blow-up technique yields a minimal nontrivial directed
cone $C$. Thus (\cite{giusti}) the dimension of the singular
set is limited by
$$H^s(\partial E\setminus\partial^\star E)=0\quad\forall s>n-7.$$
If there were another singular point apart from the origin,
the property of $C$ being a cone
would imply $H^1(\partial E\setminus\partial^\star E)=+\infty$,
contradicting $n\le7$. This shows that
there is no singular point besides the origin.
Therefore the regularity theorem implies that $C$ is a half space
and consequently we get $\partial E=\partial^\star E$.
This proof is valid for all $n\le7$. The problem in extending
this result to $\rone{}$ with $n\ge8$ is to show that the singular
set consists of at most one point. In $\R^9$ there is indeed a
counterexample which shows that there exists a nontrivial
minimal directed cone which is different from a half space.
\par
We give an alternative approach to these results that yields
more information about the structure of a directed minimal cone.
This is contained in the splitting theorem that generalizes
Corollary \ref{vierzehn}, namely, if $C$ is an open singular
minimal cone with its vertex at the origin, then $C$ is of
the form $C_0\times\R^k$, where $C_0$ is an undirected
minimal cone.

This rest of the paper is organized as follows:
In Section \ref{defkap} we introduce notations
and recall some well known facts.
In Section \ref{gerkegkap} we consider directed cones and show that
directed cones are weakly star-shaped with respect to certain
points. These properties are essential for the proof
of the regularity theorem which is stated and proved in
Section \ref{regsatzkap}. In Section \ref{splittingkap} we state and
prove the splitting theorem.
In Section \ref{anwendkap} we apply the 
regularity theorem to subgraphs, in
Section \ref{vorgegHkap} we show how it can be applied to
getting hypersurfaces
of prescribed mean curvature homeomorphic to $S^n$ for $n\le7$.
Since we use the fact that the cones arising are
directed, the results of \cite{cg} for up to seven dimensions
remain true in eight dimensions.

This paper contains unpublished results, obtained in 
1998/1999 while the author was in Heidelberg.
He would like to thank Claus Gerhardt for
many stimulating discussions about minimal cones and
their applications.

\section{Preliminaries}\label{defkap}
\begin{bez}
The characteristic function of a set $M$ will be
denoted by $\chi_M$. If we are concerned with coordinates in $\R^{n+1}$,
$(x,t)$ stands for $(x^1,\ldots,x^n,t)$ and $x=(\hat x,x^{n+1})$
abbreviates the first $n$ coordinates by $\hat x$.
Assume in the whole paper $n\ge2$.
\end{bez}

\begin{df}[almost minimal] 
A measurable set $E \subset \mathbb{R}^{n+1}$ is called almost minimal
in $\Omega\subset\R^{n+1}$, if there exists a $\lambda$ such that
$0<\lambda<1$ and for all sets $A \Subset\Omega$
there are numbers $R$ and $K$ such that
$0 < R < \dist (A, \complement\Omega)$ and $K\ge 0$ which fulfill
the inequality
$$ \int\limits_{B_\rho(x)}|D\chi_E|\le
   \int\limits_{B_\rho(x)}|D\chi_F| + K \rho^{n+2\lambda} $$
for every $x\in A$, $0 < \rho < R$ and
$E \vartriangle F \Subset B_\rho(x)$.\par
$E$ is called  (locally) minimal, if we can choose $K=0$.
\end{df}

In this paper we deal with two kinds of boundaries,
the measure-theoretical boundary and the reduced boundary,
which are defined according to \cite[p.~43]{giusti}.
Obviously these definitions are invariant, if we change
the respective sets by a set of Lebesgue measure zero.
\par
\begin{df}[(measure-theoretical) boundary]
Let $E\subset \R^{n+1}$ be a measurable set.
Define the measure-theoretical boundary
$$\partial E:=\{x\in \R^{n+1}:
  0<|B_r(x)\cap E|<|B_r(x)|\quad\forall r>0\},$$
and the measure-theoretical interior and the
measure-theoretical exterior by
$$E_\mu:=\{x\in\R^{n+1}:\exists r>0: |B_r(x)\cap E|=|B_r(x)|\}$$
and
$$\complement_\mu E:=\{x\in\R^{n+1}:\exists r>0: 0=|B_r(x)\cap E|\},$$
respectively.
\end{df}
\par
According to this definition, $E_\mu$ and $\complement_\mu E$
are open sets. $\R^{n+1}$ is the disjoint union of $E_\mu$,
$\complement_\mu E$ and $\partial E$.
Later-on we will denote $E_\mu$ resp.\ $\complement_\mu E$ simply by
$E$ resp.\ $\complement E$.
If $E$ is an almost minimal set, $E$ and $E_\mu$ lie in the
same equivalence class. This follows immediately, since
the boundary of $E$ is the union of a differentiable
manifold and of a set of $H^s$ measure zero, if $s>n-7$, and is
therefore a set of $H^{n+1}$ measure zero.
So we will always assume that an almost minimal set is open.

\begin{df}[reduced boundary]
Let $E\subset \R^{n+1}$ be a Caccioppoli set. The
reduced boundary $\partial^\star E$ consists of those points
$x\in \R^{n+1}$, for which
$$\int\limits_{B_\rho(x)}|D\chi_E|>0 \quad\forall \rho>0$$
is valid and the limit
$$\lim\limits_{\rho\downarrow 0} \nu_\rho(x)=\nu(x)$$
exists. The vector $\nu_\rho(x)$ is defined by
$$\nu_\rho(x):=\frac{\int\limits_{B_\rho(x)}D\chi_E}
  {\int\limits_{B_\rho(x)}|D\chi_E|}.$$
If $|\nu(x)|=1$, 
$\nu(x)$ is called the inner (unit) normal to $E$ at $x$.
\end{df}

\begin{df}[singularity] A point
$x\in\R^{n+1}$ is called a singularity of the measurable
set $M\subset\R^{n+1}$, if $x\in\partial M\setminus\partial^\star M$.
The set $\partial M\setminus\partial^\star M$ is called the singular
set of $M$ or the singular set of $\partial M$.
\par
A point in the boundary which is not singular is called regular.
\end{df}

\begin{df}[cone] A set
$C\subset\R^{n+1}$ is called a cone with vertex $x$, if
$$y\in C\quad\Longrightarrow\quad x+\tau(y-x)\in C\quad\forall 0<\tau.$$
Moreover, if $C\not=\emptyset$ and $C\not=\R^{n+1}$, $C$ is called
a nontrivial cone.
\end{df}

\begin{lemdef}\label{dichtgleichmass}
Let $E\subset \R^{n+1}$ be an almost minimal Caccioppoli set in the
open set $\Omega\subset\R^{n+1}$. Assume that $0\in\partial E$.
If $E_t:=\{x\in\rone{}:tx\in E\}$ converges in $L^1_{loc}(\rone{})$ to a
set $C$ for a given sequence $t_i\to0$, $t_i>0$, then $C$ is a cone
which is different from $\rone{}$ and $\emptyset$.
We define $L^1$-convergence of sets by using the corresponding
characteristic functions. 
Such a cone $C$ is called a blow-up cone of $E$ around $0$.
\end{lemdef}
\bew
Proceed as in \cite [Theorem~9.3]{giusti} and use
the estimates of \cite[p.~118]{mm} and
\cite[Proposition, p.~137]{mm}.
\qed

\begin{bem}\label{inkleinendimensionen}
In small dimensions, i.~e.\ for $n+1\le8$, the blow-up cone $C$
of an almost minimal set around a point of its boundary
has no singularity apart from the origin:
If there were another singularity, this would imply
that a half-line would be singular, because $C$ is a cone.
But this is a contradiction to the fact, that
$H^s(\partial^\star C\setminus\partial C)=0$ for $s>n-7$.
\par
Recall that the reduced boundary of a minimal set
is an analytic manifold. So, in particular, the
boundary of a minimal cone is analytic apart
from the origin, if $n+1\le8$.
\end{bem}

\begin{df}[subgraph]
Let $\varphi:A\to[-\infty,+\infty]$ be a measurable function.
Define the subgraph of the function $\varphi$ by
$${\rm sub\,} \varphi :=
\{(x, t) \in A\times \mathbb{R} : t < \varphi(x)\}.$$
\end{df}

\begin{bem}
If a cone with its vertex at the origin is a subgraph of a
function $u$, then $u$ is positive homogeneous of degree $1$.
\end{bem}

\section{Directed Cones}\label{gerkegkap}
\begin{df}[directed set]
Let $E\subset \R^{n+1}$ be a measurable set,
$y\in\R^{n+1}\setminus\{0\}$.
$$D_{\rho}(E, x):=
  \omega_{n+1}^{-1}\rho^{-n-1}\int\limits_{B_{\rho}(x)}\chi_E$$
is called approximative density of $E$ in $x$ with respect
to the radius $\rho$.
\par
A set $E$ is called directed with respect to $y$, if the maps
$$f_{x,\rho}:\R\to [0,1],\quad t\mapsto D_\rho(E,x+ty)$$
are monotone decreasing with respect to $t$ for any
$x\in\R^{n+1}$, $\rho>0$. Such $y$ is called direction of the
set $E$.
\end{df}

It is possible that a set has several linearly
independent directions.

For a measurable set $M$ we have the definition
(cf.\ e.~g.\ \cite{giusti})
\begin{align*}
x\in M&:\Longleftrightarrow \exists\rho>0:
  \int\limits_{B_\rho(x)}\chi_M=\int\limits_{B_\rho(x)}1
  \Longleftrightarrow\exists\rho>0:D_\rho(M,x)=1,\\
x\in\complement M&:\Longleftrightarrow\exists\rho>0:
  \int\limits_{B_\rho(x)}\chi_M=0
  \Longleftrightarrow\exists\rho>0:D_\rho(M,x)=0.
\end{align*}

\begin{lem}\label{untereinanderlemma}
Let $M\subset\R^{n+1}$ be a measurable set which is directed
with respect to $v$. Then for any $x\in\R^{n+1}$ and $t>0$
we have the implications
$$x\in\complement M\quad\Longrightarrow\quad x+tv\in\complement M$$
and
$$x\in M\quad\Longrightarrow\quad x-tv\in M.$$
\end{lem}
\bew
According to the definition of the
measure-theoretical complement of a set
there exists in the case $x\in\complement M$ a $\rho>0$
such that $D_\rho(M,x)=0$.
Since $M$ is directed with respect to $v$, it follows for $t>0$ that
$$D_\rho(M,x+tv)\le D_\rho(M,x)=0.$$
$D_\rho(\cdot,\cdot)$ is non-negative, so $D_\rho(M,x+tv)=0$.
This implies $x+tv\in\complement M$.
\par
In the case $x\in M$ the proof is similar:
$x\in M\Rightarrow\exists\rho>0:D_\rho(M,x)=1$.
$M$ is a directed set. Now $t>0$ implies
$D_\rho(M,x-tv)\ge D_\rho(M,x)=1$ and it follows 
that $x-tv\in M$ as above.
\qed

\begin{cor}
Let $M\subset\R^{n+1}$ be a measurable set which is directed
with respect to $v$.
If $x\in\R^{n+1}$ and $t>0$ are such that $x\in\partial M$ and
$x+tv\in\partial M$, it follows that $x+\tau v\in\partial M$ for
$0\le\tau\le t$.
\end{cor}

\begin{lem}
Let $M\subset\R^{n+1}$ be a measurable set, directed
with respect to $v$. For $x_0\in\partial M$
exactly one of the following possibilities occurs:
\begin{enumerate}
\item $\exists t\not=0:x_0+tv\in\partial M$,
\item $\forall t>0:x_0+tv\in\complement M$ and $x_0-tv\in M$.
\end{enumerate}
\end{lem}
\bew
Assume $x_0+tv\not\in\partial M$ for $t\not=0$,
i.~e.\ $x_0+tv\in M\cup\complement M$ for $t\not=0$.
Therefore we have to show that the second possibility occurs:
\par
For $t>0$ the possibility $x_0+tv\in M$ is excluded:
Assume $t>0$. Then $x_0+tv\in M$ cannot happen, because
in accordance to Lemma \ref{untereinanderlemma} $x_0+tv\in M$ implies
$x_0=(x_0+tv)-tv\in M$ contradicting $x_0\in\partial M$.
It follows $x_0+tv\in\complement M$ for $t>0$.
$x_0-tv\in\complement M$ does not occur for $t>0$ for a
similar reason. Thus the statement is proved.
\qed

\begin{lem}\label{negkomplementlemma}
Let $M\subset\R^{n+1}$ be a measurable set which is
directed with respect to $v$. Then  $\complement M$ and
$-M\equiv\{x\in\R^{n+1}:-x\in M\}$ are directed
with respect to $-v$. Especially $-\complement M$ is
again directed with respect to $v$.
\end{lem}
\bew
The statement follows at once from the equations
$$D_\rho(M,y)+D_\rho(\complement M,y)=1$$
and
$$D_\rho(-M, -x-\tau v)=D_\rho(M, x+\tau v)$$
which are valid for any $x\in\rone{}$ and $\rho>0$.
\qed

\begin{lem}\label{cistgerichteterkegel}
Let $E\subset\R^{n+1}$ be representable as a subgraph.
If $t_i\downarrow0$ is as before such that
$E_{t_i}\equiv\{x\in\R^{n+1}:t_ix\in E\}$ converges in
$L^1_{loc}(\R^{n+1})$ to a cone $C$, then $C$ is
directed with respect to $e_{n+1}$.
\end{lem}
\bew
As a subgraph, $E$ is directed with respect to $e_{n+1}$
(cf.\ Remark \ref{subgleichger}).
If $C$ were not directed with respect to $e_{n+1}$, there
would be $x\in\R^{n+1}$, $\rho>0$ and $t>0$ such that
$$\int\limits_{B_\rho(x)}\chi_C<\int\limits_{B_\rho(x+te_{n+1})}\chi_C.$$
In the same way as for $E$ we get that $E_{t_i}$ is directed with
respect to $e_{n+1}$. Because of the convergence
$E_{t_i}\to C$ in $L^1_{loc}(\R^{n+1})$ we immediately get
a contradiction to the inequality above.
\qed

\begin{df}[weakly star-shaped]
A set $M\subset\R^{n+1}$ is called weakly star-shaped with
center $x$, if for all $z\in M$ we have
$x+\tau(z-x)\in M$ if $0<\tau\le 1$.
\end{df}

The following Lemma is essential for the proof of the
regularity theorem:
\begin{lem}\label{sternlemma}
Let $C$ be an open cone with vertex at the origin which is
directed with respect to $e_{n+1}$.
Then $C$ is weakly star-shaped with center $x=(0,-t)$ for
all $t>0$.
\end{lem}
\bew
Let $y=(\hat y,y^{n+1})\in C$ and $\tau$ with $0<\tau<1$ be arbitrary.
We show that $x+\tau(y-x)=(\tau\hat y,-t+\tau(y^{n+1}+t))$ is an
element of $C$. Being a cone, $C$ contains $(\tau\hat y,\tau y^{n+1})$
because of $y\in C$. $-t(1-\tau)$ is negative, so Lemma
\ref{untereinanderlemma} implies
$(\tau\hat y,\tau y^{n+1})-t(1-\tau)(0,1)=
  (\tau\hat y,-t+\tau(y^{n+1}+t))\in C$.
\qed

\begin{bem}\label{subgleichger}
For a measurable set $M\subset\R^{n+1}$ such that
$H^{n+1}(\partial M)=0$ the following two statements
are equivalent:
\begin{enumerate}
\item $M$ is the subgraph of a measurable function $u$,
\item $M$ is directed with respect to $e_{n+1}$.
\end{enumerate}
The measurable function $u$ in (i) is given by
$$u(\hat x):=\sup\{t\in\R:(\hat x,t)\in M\}.$$
\end{bem}
\bew
``(i)$\Longrightarrow$(ii)'':\\
Let $M={\rm sub\,}u$ be a given set.
It follows that $\chi_M(x,t)=1$ for $t<u(x)$ and
$\chi_M(x,t)=0$ for $t>u(x)$, i.~e.\ $\chi_M(x,t)$ is
monotone decreasing with respect to $t$ for any $x\in\R^n$.
By integrating, we get for any $\rho>0$,
$\tau>0$ and $x\in\R^{n+1}$
$$\int\limits_{B_\rho(x)}\chi_M(z)dz
\ge\int\limits_{B_\rho(x)}\chi_M(z+\tau e_{n+1})dz
=\int\limits_{B_\rho(x+\tau e_{n+1})}\chi_M(y)dy.$$
Therefore $M$ is a directed set.\par
``(ii)$\Longrightarrow$(i)'':\\
Define for $\hat x\in\R^n$
$$u(\hat x):=\sup\left\{t\in\R:(\hat x,t)\in M\right\}.$$
We remark that $(\hat x,t)\in M$ is equivalent to the
existence of a $\rho>0$ such that $D_\rho(M,(\hat x, t))=1$.
Since $M$ is a measurable set the function $u$ is measurable.
Define $U:={\rm sub\,}u$.\par
{\bf a)} $M\subset U$:\\
Assume that $(\hat x,t)\in M$. Since the supremum in the definition
of $u(\hat x)$ is not assumed, it follows $u(\hat x)>t$ and
therefore $(\hat x,t)\in U$.\par
{\bf b)} $U\subset M$:\\
Assume $(\hat x,t)\in U$, i.~e.\ $u(\hat x)>t$.
The definition of $u(\hat x)$ implies that there exists a
$\tau$ such that $u(\hat x)>\tau>t$ and $(\hat x,\tau)\in M$.
Now Lemma \ref{untereinanderlemma} yields $(\hat x,t)\in M$, since
$t-\tau<0$ and $(\hat x,\tau)+(t-\tau)(0,1)=(\hat x,t)$.
\qed

\begin{df}[cone of directions]\label{defconeofdirections}
Let $E\subset\rone{}$ be a measurable set. $\Dir(E)$, the cone
of directions, is defined to be the set of all directions of
$E$ together with the origin.
\end{df}

\begin{rem}
Definition \ref{defconeofdirections} is equivalent to
\begin{align*}
\Dir(E):=&\{y\in\rone{}:f_{x,\rho}:\R\to[0,1],
  t\mapsto D_\rho(E,x+ty)\\
  & \text{ is monotone decreasing for any }x\in\rone{}
    \text{ and any }\rho>0\}.\\
\end{align*}
\end{rem}

\begin{lem}\label{kegelistmonoid}
Let $E\subset\rone{}$ be a measurable set. Then $\Dir(E)$ is a closed
cone which is also closed under addition.
\end{lem}
\bew\par
{\bf (i)} The fact, that $\Dir(E)$ is a cone follows immediately
from the Definition \ref{defconeofdirections}.\par
{\bf (ii)} $\Dir(E)$ is closed under addition:\\
Let $y_1,y_2\in\Dir(E)$ be arbitrary. According to the definition this
is equivalent to
$$D_{\rho_1}(E,x_1+t_1y_1)\le D_{\rho_1}(E,x_1+\tau_1y_1)$$
and
$$D_{\rho_2}(E,x_2+t_2y_2)\le D_{\rho_2}(E,x_2+\tau_2y_2)$$
for any $x_1,x_2\in\rone{}$, $\rho_1,\rho_2>0$ and
$t_1,t_2,\tau_1,\tau_2\in\R$ such that $t_1\ge\tau_1$ and
$t_2\ge\tau_2$.
We have to show now, that
$$D_\rho(E,x+t(y_1+y_2))\le D_\rho(E,x+\tau(y_1+y_2))$$
for any $x\in\rone{}$, $\rho>0$ and $t,\tau\in\R$ such that $t\ge\tau$.
We choose now $x_1=x+ty_2$, $t_1=t$, $\rho_1=\rho$, $\tau_1=\tau$,
$x_2=x+\tau y_1$, $t_2=t$, $\rho_2=\rho$, $\tau_2=\tau$
and deduce from the inequalities above
\begin{align*}
D_\rho(E,(x+ty_2)+ty_1)\le&D_\rho(E,(x+ty_2)+\tau y_1)\\
=&D_\rho(E,(x+\tau y_1)+t y_2)\\
\le&D_\rho(E,(x+\tau y_1)+\tau y_2)\\
\end{align*}
verifying the claimed inequality.\par
{\bf (iii)} $\Dir(E)$ is a closed set:\\
Let $y_i\in\Dir(E)$ for $i\in\N$ such that $y_i\to y$ as
$i\to\infty$. We have to show $y\in\Dir(E)$.
Assume in contrast $y\not\in\Dir(E)$. Then there exists $x\in\rone{}$,
$\rho>0$, $t,\tau\in\R$ such that $t\ge\tau$ and
$$D_\rho(E,x+ty)-D_\rho(E,x+\tau y)>0.$$
\OE we can assume $\tau=0$. $y_i\in\Dir(E)$ implies
$$D_\rho(E,x+ty_i)-D_\rho(E,x)\le0.$$
Since $D_\rho(E,x+ty_i)$ converges to $D_\rho(E,x+ty)$
as $i$ tends to infinity we get a contradiction and the statement
follows.
\qed

\section{Regularity Theorem}\label{regsatzkap}

\begin{df}\label{aufeinerseite}
An open set $M\subset\R^{n+1}$ is said to lie on one side
of a hyperplane $T$ if an adjusted rotation and translation of
the coordinate system yields the following situation
$$x=(\hat{x},x^{n+1})\in T\quad\Longleftrightarrow\quad x^{n+1}=0,$$
$$x=(\hat{x},x^{n+1})\in M\quad\Longrightarrow\quad x^{n+1}>0.$$
\end{df}

\begin{lem}\label{c1kegelisthalbraum}
Let $C\subset\R^{n+1}$ be a cone which is representable as a
subgraph of a $C^1$-function $u$ in a neighborhood of its vertex.
Then  $\partial C$
and the tangential hyperplane $T$ to $\partial C$ at
the vertex of $C$ coincide and $C$ is a half space.
\end{lem}
\bew
By a translation we can assume \oe that the vertex of the cone
is the origin. Since $C$ is a cone, $u$ is positive
homogeneous of degree $1$ and has a well-defined extension
(\oe $u$)
of the same homogenity which is defined on the whole $\R^n$.
The subgraph of $u$ is $C$.
Let $v\in\R^n$ be arbitrary. It follows
$$\langle Du(0),v\rangle=\lim\limits_{t\to0}\frac{u(tv)-u(0)}{t}.$$
If we take into account $u(0)=0$ and use the fact that $u$ is a
positive homogeneous function, we deduce for $t>0$
$$\langle Du(0),v\rangle=
  \lim\limits_{t\to0}\frac{t\cdot u(v)}{t}=u(v).$$
The left-hand side of this equality is linear with respect to $v$.
Hence $u$ is a linear function. Thus $C$ is a subgraph of a
linear function and the statement follows.
\qed

\begin{lem}\label{zweipktebene}
Let $C\subset\R^{n+1}$ be an open minimal cone with its vertex
at the origin. If $C$ is directed with respect to $e_{n+1}$
and there exists a $t>0$ such that $z:=(0,t)\in\complement C$ and
$(0,-t)\in C$, then $\partial C$ is a hyperplane.
\end{lem}
\bew
\OE we assume $n\ge7$.
Choose $r>0$ such that $D_r(C,z)=0$ and $D_r(C,-z)=1$.
Define a cone $K$ by
$$K:=\left\{(\hat x,-\tau):\tau>0, \tau r>|\hat x|t\right\}.$$\par
It follows\par
{\bf (i)} $K\subset C$ and $-K\subset\complement C$:\\
Let $(\hat x,-\tau)\in K$ be arbitrary.
Since $K$ is a cone, we get $\left(\hat x\frac{t}{\tau},
-\tau\frac{t}{\tau}\right)=\left(\hat x\frac{t}{\tau},-t\right)\in K$.
Now $\left|\hat x\frac{t}{\tau}\right|<r$ implies
$\left(\hat x\frac{t}{\tau},-t\right)\in B_r(-z)$.
Choose $s>0$ such that
$B_s\left(\left(\hat x\frac{t}{\tau},-t\right)\right)\subset
B_r(-z)$. The equality $D_r(C,-z)=1$ implies
$D_s\left(C,\left(\hat x\frac{t}{\tau},-t\right)\right)=1$
and therefore we get
$\left(\hat x\frac{t}{\tau},-t\right)\in C$.
$C$ is a cone, so we deduce $(\hat x,-\tau)\in C$.
As $(\hat x,-\tau)\in K$ was arbitrary, so we get $K\subset C$.
\par
In the same way
$-K:=\left\{(\hat x,\tau):\tau>0,\tau r>|\hat x|t\right\}
\subset\complement C$ is proved.\par
{\bf (ii)} Representability of $\partial C$ as a graph:\\
We will show, that $\partial C$ can be represented as a
subgraph over $\R^n\setminus\Sigma$ besides a set which has
$H^{n-5}$-measure zero. $\Sigma$ is a closed set of
$H^{n-6}$-measure zero.
\par
For $s>n-7$ we get $H^s(\partial C\setminus\partial^\star C)=0$
by Theorem \cite[Theorem~11.8, p.~134]{giusti}.
This implies especially $H^{n-6}(\partial C\setminus\partial^\star C)=0$.
Define $\pi:\R^n\times\R\to\R^n$ by $\pi((\hat x,x^{n+1}))=\hat x$
and $\Sigma:=\pi(\partial C\setminus\partial^\star C)$.
Since $\pi$ is Lipschitz continuous we deduce
$H^{n-6}(\Sigma)=0$, hence
$H^{n-5}(\Sigma\times\R)=0$ and
$H^{n-5}((\Sigma\times\R)\cap\partial C)=0$.
$\Sigma$ is a closed set, because $\pi$ is continuous and
in view of (i) $(\partial C\setminus\partial^\star C)\cap
(F\times\R)$ is compact for any compact set $F\subset\R^n$.
\par
Let $\hat x\in\R^n\setminus\Sigma$ be arbitrary. As $K\subset C$,
$\tau>\frac{|\hat x|t}{r}$ implies
$(\hat x,\tau)\in\complement C$ and
$(\hat x,-\tau)\in C$.
$\{\hat x\}\times\R$ is connected, $C$, $\complement C$ und
$\partial C$ are disjoint and we have
$C\cup\complement C\cup\partial C=\R^n$.
Regarding the fact, that $C$ and $\complement C$
are open sets we deduce
$(\{\hat x\}\times\R)\cap\partial C\not=\emptyset$,
because a connected set is not
the disjoint union of two non-empty open sets, in this case
$(\{\hat x\}\times\R)\cap C$ with $(\hat x,-\tau)\in C$ and
$(\{\hat x\}\times\R)\cap\complement C$ with
$(\hat x,\tau)\in\complement C$ for all
$\tau>\frac{|\hat x|t}{r}$. The set
$(\{\hat x\}\times\R)\cap\partial C$ is bounded because
$(\hat x,\tau)\in\complement C$ and
$(\hat x,-\tau)\in C$ for all $\tau>\frac{|\hat x|t}{r}$.
$C$ is directed with respect to $e_{n+1}$ and $\partial C$ is closed,
so Lemma \ref{untereinanderlemma} implies
$(\{\hat x\}\times\R)\cap\partial C=\{\hat x\}\times I$ for a
compact interval $I$. The boundary of $C$ is analytic in the
complement of $\Sigma\times\R$ and so $I$ consists of exactly
one point. Therefore we have a function
$u\in L^1_{loc}(\R^n\setminus\Sigma)$, whose graph
coincides with $\partial C\setminus(\Sigma\times\R)$ and we
have the equality $C\cap((\R^n\setminus\Sigma)\times\R)={\rm sub\,} u$.
Observe, however, that $u$ is not automatically analytic, because
the boundary $\partial C$ in the complement of $\Sigma\times\R$
is only analytic as a manifold.\par
{\bf (iii)} Regularity of  $u$:\\
Let $\Omega\Subset\R^n\setminus\Sigma$ be an open ball.
Since the subgraph of $u|_\Omega$ has finite perimeter in
$\Omega\times\R$, i.~e.\
$\int\limits_{\Omega\times\R}|D\chi_{\sub u}|<\infty$,
and $u\in L^\infty(\Omega)$, we can use
\cite[Theorem~1, p.~317]{gia1} to deduce $u\in BV(\Omega)$.
According to \cite{cg1} we get $u\in C^{0,1}(\Omega)$.
But $u$ is also a weak solution of the minimal surface equation
and this implies $u\in C^2(\Omega)$.
$\Omega\Subset\R^n\setminus\Sigma$ was an arbitrary open ball.
Thus we deduce $u\in C^2(\R^n\setminus\Sigma)$. 
We have $H^{n-6}(\Sigma)=0$, so Theorem \cite[Theorem 16.9]{giusti}
can be applied, i.~e.\ $u$ can be extended to a
function in $C^2(\R^n)$ solving the minimal
surface equation. $\sub u$ is an element of the
$L^1$-equivalence class of $C$ independent of the choice
of the extension, because another extension changes $\sub u$
at most by a subset of $\Sigma\times\R$ and
$H^{n-5}(\Sigma\times\R)=0$.
\par
{\bf (iv)} $\partial C$ is a hyperplane:\\
This follows immediately from Lemma \ref{c1kegelisthalbraum}.
\qed

The following Lemma is - besides Lemma \ref{sternlemma} -
the essential part of the proof of the regularity theorem:
\begin{lem}\label{einseitenlemma}
Let $M\subset\R^{n+1}$ be an open 
weakly star-shaped set with center $x_0\in\partial M$.
If the boundary of $M$ is of class $C^1$ in a
neighborhood of $x_0$, then $M$ lies on one side
of the tangential hyperplane $T=T_{x_0}\partial M\subset\R^{n+1}$.
\end{lem}
\bew
By a translation of the coordinates we can assume that
$x_0=(0,0)\in\R^n\times\R$. Since $\partial M$ is
of class $C^1$ in a neighborhood of $x_0$, we reach the following
situation by a suitable rotation of the coordinates:
There exist $R>0$ and $u\in C^1({\hat{B}}_{2R})\equiv
C^1(\{\hat{x}\in\R^n:|\hat{x}|<2R\})$ such that
$u(0)=0$, $Du(0)=0$ and for any $x=(\hat{x},x^{n+1})\in
{\hat{B}}_R\times(-R,R)$ we have the three equivalences
$$x\in\partial M\quad\Longleftrightarrow\quad x^{n+1}=u(\hat{x}),$$
$$x\in M\quad\Longleftrightarrow\quad x^{n+1}>u(\hat{x})$$ and
$$x\in\complement M\quad\Longleftrightarrow\quad x^{n+1}<u(\hat{x}).$$
By construction we have for any $x\in\R^{n+1}$
$$x\in T\quad\Longleftrightarrow\quad x^{n+1}=0.$$
\par
Assume now, the assertion would be false, more precisely: The chosen
coordinates were not fulfilling the conditions of
Definition \ref{aufeinerseite}. That is, the condition
$$x=(\hat{x},x^{n+1})\in M\quad\Longrightarrow\quad x^{n+1}>0$$
is violated for an $x\in\R^{n+1}$. As $M$ is an open set, 
there exists a point
$y_0=({\hat{y}}_0, y^{n+1}_0)\in M$ with $y^{n+1}_0<0$.
First, we observe that ${\hat{y}}_0=0$ cannot occur,
for $M$ is weakly star-shaped with center $x_0=(0,0)$ and
we have for $y_0\in M$ that $\tau y_0$ is in $M$
for $0<\tau\le1$. This implies for $|\tau y_0|<R$, i.~e.\ for
sufficiently small $\tau>0$, we get the inequality
$$\tau y_0^{n+1}<0=u(0)=u({\hat{y}}_0)$$
which contradicts the assumption
$$x\in M\quad\Longleftrightarrow\quad x^{n+1}>u(\hat{x})$$
for $x\in{\hat{B}}_R\times(-R,R)$.
Assume therefore  ${\hat{y}}_0\not=0$.
Using again the fact that $M$ is weakly star-shaped, we get
$\tau y_0=(\tau{\hat{y}}_0,\tau y^{n+1}_0)\in M$ for $0<\tau\le1$.
Choose now $\varepsilon>0$ such that
$$0<\tau<\varepsilon\quad\Longrightarrow
  \quad\tau y_0\in{\hat{B}}_R\times(-R,R).$$
In ${\hat{B}}_R\times(-R,R)$, the set $M$ was characterized by
the equivalence
$$x=(\hat{x},x^{n+1})\in M\quad\Longleftrightarrow
  \quad x^{n+1}>u(\hat{x}).$$
This implies $u(\tau{\hat{y}}_0)<\tau y_0^{n+1}$
for $0<\tau<\varepsilon$.
If we take into account that $u(0)=0$, $|{\hat{y}}_0|\not=0$
and $y_0^{n+1}<0$, we deduce
$$u(\tau{\hat{y}}_0)<\tau y_0^{n+1}<0,$$
$$|u(\tau{\hat{y}}_0)-u(0)|=-u(\tau{\hat{y}}_0)>-\tau y_0^{n+1},$$
and finally
$$\left|\frac{u(\tau{\hat{y}}_0)-u(0)}{\tau|{\hat{y}}_0|}\right|
  >\frac{-y_0^{n+1}}{|{\hat{y}}_0|}$$
for all $0<\tau<\varepsilon$.
The right-hand side of this inequality is independent of $\tau$
and positive, so we get a contradiction for $\tau\to0$
to the fact that the coordinates were chosen such that
$Du(0)=0$.
\qed

\begin{satz}[Regularity Theorem]\label{regsatz}
Let $C\subset\R^{n+1}$ be an open nontrivial minimal cone which
is directed with respect to $e_{n+1}$. Assume that the vertex of
$C$ is the origin and that $C$ has at most one singularity,
i.~e.\ $C$ is at most singular in the origin.
Then $\partial C$ is a hyperplane and $C$ is a half space.
\end{satz}
For $n\le6$ there is nothing to be proved, the statement is true
according to \cite[Theorem~10.11, p.~127]{giusti} even for non-directed
cones.\par
\bew
$\partial C$ is regular apart from the origin, so
$\partial C\setminus\{0\}$ is analytic.
\par
In view of Lemma \ref{zweipktebene}
and due to the fact that $C$ is directed, it suffices to consider the
case that $x_0=(0,t)\in\partial C$ for some $t\not=0$.\par
We may assume that $t<0$. Otherwise consider
$-\complement C\equiv\{x\in\R^{n+1}:-x\in\complement C\}$
instead of $C$. According to Lemma \ref{negkomplementlemma},
this set is directed with respect to $e_{n+1}$.
\par
Now Lemma \ref{sternlemma} implies that $C$ is weakly
star-shaped with center $x_0$. $\partial C$ is regular
apart from the origin. Therefore we have a well-defined
tangential hyperplane $T$ at $\partial C$. With the help
of Lemma \ref{einseitenlemma} we conclude, that $C$ lies
on one side of $T$. Since $C$ is nontrivial,
we get $0\in\partial C$.
Finally \cite[Theorem~15.5, p.~174]{giusti}
implies that $C$ is a half space and $\partial C=T$
is a hyperplane.
\qed

\begin{cor}\label{regularitaetssatzcorollar}
Let $C\subset\R^{n+1}$ be an open nontrivial minimal cone which is
directed with respect to $e_{n+1}$.
If $n+1\le8$ then $\partial C$ is a hyperplane and $C$
is a half space.
\end{cor}
\bew
Remark \ref{inkleinendimensionen} guarantees that $C$
has at most one singularity which is at the vertex if it exists.
Thus our Regularity Theorem \ref{regsatz} yields the statements.
\qed

\begin{bem}
In the Regularity Theorem \ref{regsatz},
we assumed that $C$ is regular outside its vertex.
In the proof, however, we
use only the fact, that there is a $t>0$ which satisfies the
following two conditions:\\
\quada {\rm (i)} $(0,t)\in\partial C$ $\Longrightarrow$
  $\partial C$ is a $C^1$-manifold in a neighborhood of $(0,t)$.\\
\quada {\rm (ii)} $(0,-t)\in\partial C$ $\Longrightarrow$
  $\partial C$ is a $C^1$-manifold in a neighborhood of $(0,-t)$.\\
Therefore we get the following corollary:
\end{bem}

\begin{cor}\label{minregsatzcorollar}
Let $C\subset\rone{}$ be an open nontrivial minimal cone which is
directed with respect to $e_{n+1}$ and which satisfies the two
conditions (i) and (ii) stated above for a positive $t>0$.
(If one of the points $(0,t)$ or $(0,-t)$ is an element of
$\partial C$, we need in fact only the respective condition.)
Then $C$ is a half space.
\end{cor}

\begin{bem}
Conditions (i) and (ii) in Corollary \ref{minregsatzcorollar}
can alternatively be replaced by one of the following conditions:\\
\quada {\rm (i)} There is a direction $w$ of $C$ which fulfills
  $w\not\in\Sigma$ and $-w\not\in\Sigma$ or
  $w\in\partial C\setminus\Sigma$ or
  $-w\in\partial C\setminus\Sigma$.\\
\quada {\rm (ii)} There is one direction $w$ of $C$
  which fulfills $w\not\in\langle\Sigma\rangle$.\\
\quada {\rm (iii)} $\dim\langle D\rangle>\dim\langle\Sigma\rangle$.\\
Here $D$ is the set of all directions and $\Sigma$ is the singular
set of $C$. The brackets denote the vectorspace which is
spanned by the respective set.
\end{bem}
\bew
This follows at once from Corollary \ref{minregsatzcorollar},
if we use the fact, that $\partial C\setminus\Sigma$ is an
open $C^1$-manifold.
\qed

\begin{cor}\label{vierzehn}
Let $C\subset\rone{}$ be an open nontrivial minimal cone with
vertex at the origin. Suppose $C$ has $k$ linearely
independent directions, $k\in\N$. Then $C$ is a half-space
provided that $n+1-k\le7$ or
$H^{k}(\partial C\setminus\partial^\star C)=0$.
\end{cor}
\bew
In the first case ($n+1-k\le7$) we get
$H^k(\partial C\setminus\partial^\star C)=0$, for $k>n-7$
and $C$ is an almost minimal set, so the second case includes
the first one.
\par
Since $C$ has $k$ linearely independent directions we
deduce according to Lemma \ref{kegelistmonoid} that
$H^k(\Dir(C))=\infty$.
But we assumed that the $k$-dimensional Hausdorff-measure
of the singular set of $C$ vanishes. So we can find
$x\in\reu$ such that
$x\not\in\partial C\setminus\partial^\star C$,
$-x\not\in\partial C\setminus\partial^\star C$ and
$x\in\Dir(C)$. Then we apply
Corollary \ref{minregsatzcorollar}.
\qed

\section{Applications}\label{anwendkap}

\begin{satz}[regularity for subgraphs]\label{subgraphanwend}
Assume that $E\subset\R^{n+1}$, $n\le 7$, is a subgraph and
almost minimal (with constant $\lambda$) in $\Omega$,
$\Omega\subset\R^{n+1}$ open,
then we get $\partial E\cap\Omega=\partial^\star E\cap\Omega$
and $\partial E\cap\Omega$ is a $C^{1,\lambda}$-manifold.
\end{satz}
\bew
According to the definitions we get
$\partial^\star E\subset\partial E$.
Let $x_0\in\partial E\cap\Omega$ be an arbitrary point.
Show $x_0\in\partial^\star E$:
By virtue of a translation we can assume that $x_0=0$.
Then \cite[Proposition, p.~137]{mm}
guarantees that there is a sequence $t_i\to0$, $t_i>0$, such that
$E_{t_i}:=\{x\in\R^{n+1}:t_ix\in E\}$
converges for $i\to\infty$ in $L^1_{loc}(\R^{n+1})$
to a minimal cone $C$.
Lemma \ref{dichtgleichmass} ensures that $C$
is nontrivial. Then we get, according to the
quoted proposition
$$0\in\partial^\star C\quad\Longleftrightarrow\quad
  x_0=0\in\partial^\star E.$$
Assume that $C$ is open. We know 
\cite[Lemma~16.3, p.~184]{giusti}
that $C$ is representable as a subgraph.
Thus Theorem \ref{regsatz} implies that $\partial C$ is a
hyperplane, and therefore we get $\partial C=\partial^\star C$.
Using the equivalence from above, we get
$x_0\in\partial^\star E$. Finally, \cite{tam2,tamanini}
yields that $\partial E\cap\Omega$ is a
$C^{1,\lambda}$-manifold. The theorem is proved.
\qed

\begin{thm}\label{simonthm}
Simons' cone
$$K^{2m}:=\left\{x\in\R^{2m}:\sum\limits_{i=1}^m\left(x^i\right)^2<
\sum\limits_{i=m+1}^{2m}\left(x^i\right)^2\right\}$$
is minimal for $m\ge4$.
\end{thm}
\bew
\cite[Theorem~16.4, p.~185]{giusti}
\qed

\begin{lem}\label{neundimensionalesgegenbeispiel}
Theorem \ref{subgraphanwend} is false in $\R^9$, i.~e.\ if
we replace the assumption $E\subset\rone{}$, $n\le7$, by
$E\subset\rone{}$, $n\le8$, we cannot prove
$\partial E\cap\Omega=\partial^\star E\cap\Omega$.
\end{lem}
\bew
Let $K:=K^8$ be the minimal cone defined in Theorem \ref{simonthm}.
$K\times\R\subset\R^9$ is obviously the subgraph of the function $u$
defined by
$$u(\hat x):=\left\{\begin{array}{r@{\quad:\quad}l}
  +\infty & \hat x\in K\\
  -\infty & \hat x\not\in K\\
  \end{array}\right.$$
for all $\hat x\in\R^8$.
The set $K\times\R$ is minimal \cite[Example~16.2, p.~183]{giusti}.
In this example, the measure-theoretical and the reduced boundary
differ, we have $\partial(K\times\R)=
\partial^\star(K\times\R)\dot{\cup}(\{0\}\times\R)$.
This equality follows,
because the reduced boundary of an almost minimal
set $E\subset\rone{}$ consists exactly of those points $x$ of the
boundary for which $\partial E$ is a $C^1$-manifold in
a neighborhood of $x$.
\qed

\def\ds{{d\overline{s}}}
\def\tN{{\tilde{N}}}
\def\gab{{\overline{g}_{\alpha\beta}}}
\def\g{{\overline{g}}}
\def\gh{{\hat g}}

\section{Prescribed Mean Curvature}\label{vorgegHkap}
The following problem is considered in \cite{cg}:\\
In a complete locally conformally flat $(n+1)$-dimensional
Riemannian manifold $N$ we look for a closed
hypersurface $M$ which is homeomorphic to $S^n$ and
has prescribed mean curvature $f$, $f\in C^{0,1}(N)$, i.~e.\
the equation $H|_M=f(x)$ shall be solved by a hypersurface
of class $C^{2,\alpha}$. The hypersurface is looked for in
an open connected relatively compact subset
$\Omega$ of $N$, which is also regarded - using a diffeomorphism -
as a subset of $\rone{}$. The boundary of $\Omega$ consists of
two components $M_1$ and $M_2$, which are given
in Euclidean polar coordinates $(x^\alpha)_{0\le\alpha\le n}$,
$x^0\equiv r$, $u_i\in C^{2,\alpha}(S^n,(0,\infty))$,
as graphs over $S^n$:
$M_i={\rm graph\,}\left.u_i\right|_{S^n}=
  \{(z,u_i(z)):z\in S^n\}$.
$M_1$ and $M_2$ act as barriers, i.~e.\ they satisfy
$H|_{M_1}\le f$ and $H|_{M_2}\ge f$, where the respective
unit normal vector $(\nu^\alpha)$ is chosen such that the component
$\nu^0$ is negative.
\par
In the cited paper it is proven, that under the 
assumptions stated above,
such a hypersurface $M$ exists provided that $n\le6$.
In this chapter we extend the proof and show that such
a hypersurface $M$ exists up to $n=7$.

\begin{thm}\label{hyperflaechenbisdimensionacht}
The problem
``Find a closed hypersurface $M\subset\overline{\Omega}$
of class $C^{2,\alpha}$ such that $\left.H\right|_M=f$, which
is homeomorphic to $S^n$.'' has a solution
if $n\le7$.
\end{thm}
We need some Lemmata:

\begin{lem}\label{randinvarianzlemma}
The metric product $S^n\times\R$ is locally conformally equivalent
to $\rone\setminus\{0\}$.
\end{lem}
\bew
Let $\tN=S^n\times\R$ be the metric product of $S^n$ and $\R$.
The metric of  $\tN$ is given by
$$\ds_\tN^2=dt^2+\sigma_{ij}dx^idx^j,$$
where $(x^i)_{1\le i\le n}$ are coordinates of $S^n$ and $t\in\R$.
Now we identify the manifold $\tN$ with its image in $\rone{}$ under the
diffeomorphism
\begin{align*}
\Psi:S^n\times\R\to& \reu,\\
(x,t)\mapsto& (x,e^t)\equiv(x,r).\\
\end{align*}
$\Psi$ is an isometry, if we equip $\reu$ with the metric
$$\frac{1}{r^2}dr^2+\sigma_{ij}dx^idx^j.$$
$((x^i)_{1\le i\le n}, r)$ are polar coordinates of $\reu$.
Now we assume that
$(\tN,\ds_\tN^2)=(\reu,\frac{1}{r^2}dr^2+\sigma_{ij}dx^idx^j)$
and deduce
\begin{align*}
\ds^2_\tN
=& \frac{1}{r^2}dr^2+\sigma_{ij}dx^idx^j\\
=& \frac{1}{r^2}(dr^2+r^2\sigma_{ij}dx^idx^j)\\
=& \frac{1}{r^2}\ds^2_\rone{}.
\end{align*}
In the last equality we use the fact that $dr^2+r^2\sigma_{ij}dx^idx^j$
is a representation of the standard metric of $\reu$ in
polar coordinates. This equation yields that
$(\tN,\ds_\tN^2)$ is a locally conformally flat Riemannian manifold.
\qed

\begin{lem}\label{alssubgraphdarstellbarlemma}
Let $E\subset\rone{}$ be an almost minimal set which is representable
as a subgraph over $S^n$, i.~e.\
$$E=\{(\hat x,t)\in S^n\times\R^+\subset\reu:t<u(\hat x)\}$$
is a subgraph in polar coordinates.
Assume $u\ge c>0$.
Let $z_0=(0, z_0^{n+1})\subset\rone{}$ be an arbitrary point
such that $z_0^{n+1}>0$.
If $E_{t_i}=\{y\in\rone{}:z_0+t_i(y-z_0)\in E\}$
converges for a given sequence $t_i\to0$, $t_i>0$,
$i\in\N\setminus\{0\}$, in $L^1_{loc}(\rone{})$
to a cone $C$, then $C$ is representable as a
subgraph over $\R^n\equiv{\langle z_0\rangle}^\perp$.
\end{lem}
\bew
Identify $\R^n$ with $\langle z_0\rangle^\perp$ and introduce
orthogonal coordinates.
Show that there is no $\hat x\in\R^n$ such that
$(\hat x, t)\in C\cup\partial C$ and $(\hat x, \tau)\in\complement C$
with $t>\tau$ or $(\hat x, t)\in C$ and $(\hat x, \tau)\in\partial C$
with $t>\tau$.
\par
{\bf (i)} Assume that $z_0^{n+1}=1$. In order to get an
easier representation in coordinates, we translate such that $z_0=0$.
Then we have
$$E_{t_i}=\{y\in\rone{}:t_iy\in E\}.$$
Now, $E$ is representable as a subgraph over $S^n$ with center
$(0,-1)$ and $E_{t_i}$ is representable as a subgraph
over $S^n$ with center $(0,-\frac{1}{t_i})$.
\par
{\bf (ii)} $(\hat x, t)\in C$ and
$(\hat x,\tau)\in\complement C$ $\Longrightarrow$ $t<\tau$:\\
Assume that there is a $\hat x_0\in\R^n$, such that
$(\hat x_0,t_0)\in C\cup\partial C$,
$(\hat x_0,\tau_0)\in\complement C$ and $t_0>\tau_0$.
Then there exists $\rho>0$ such that
$D_\rho(C,(\hat x_0,\tau_0))=0$. On the other hand, we have
$D_\rho(C,(\hat x_0,t_0))>2\varepsilon>0$.
Since $E_{t_i}\to C$ in $L^1_{loc}(\rone{})$, we get the inequality
$D_\rho(E_{t_i},(\hat x_0,t_0))>\varepsilon>0$ for sufficiently
large $i$.
Since each $E_{t_i}$ is representable as a subgraph over $S^n$ with
center $\left(0,-\frac{1}{t_i}\right)$, we deduce that
$$(0,\infty)\ni s\mapsto
  D_{s\rho}\left(E_{t_i},\left(0,-\frac{1}{t_i}\right)+sx\right)$$
is monotone decreasing for any $x\in\rone{}$,
because
\begin{align*}
&\chi_{E_{t_i}}\left(\left(0,-\frac{1}{t_i}\right)+s_0y\right)=0,\,
  s_0>0,\,y\in\rone{}\\
\Longrightarrow&
  \chi_{E_{t_i}}\left(\left(0,-\frac{1}{t_i}\right)+sy\right)=0
  \quad\forall s\ge s_0\\
\end{align*}
implies
$$s_0^{-n-1}\int\limits_{B_{s_0\rho}
  \left(\left(0,-\frac{1}{t_i}\right)+s_0x\right)}\chi_{E_{t_i}}\ge
  s^{-n-1}\int\limits_{B_{s\rho}
  \left(\left(0,-\frac{1}{t_i}\right)+sx\right)}
  \chi_{E_{t_i}}\,\forall s\ge s_0.$$
Now, for large $i$ the number
$s_i:=\frac{\tau_0+\frac{1}{t_i}}{t_0+\frac{1}{t_i}}$
satisfies $0<s_i<1$, and therefore we get for such $i$
\begin{align*}
0<\varepsilon<&D_\rho(E_{t_i},(\hat x_0,t_0))\\
=&D_{1\rho}\left(E_{t_i},\left(0,-\frac{1}{t_i}\right)
  +1\left(\hat x_0,t_0+\frac{1}{t_i}\right)\right)\\
\le&D_{s_i\rho}\left(E_{t_i},\left(0,-\frac{1}{t_i}\right)+
  s_i\left(\hat x_0,t_0+\frac{1}{t_i}\right)\right)\\
=&D_{s_i\rho}(E_{t_i},(s_i\hat x_0,\tau_0))
\end{align*}
On the other hand, $s_i$ converges to $1$ as $i$ tends
to infinity, so we get
$D_\rho(C,(\hat x_0,\tau_0))\ge\varepsilon>0$;
this inequality contradicts $D_\rho(C,(\hat x_0,\tau_0))=0$.
\par
{\bf (iii)} $(\hat x,t)\in C)$ and $(\hat x,\tau)\in\partial C$
$\Longrightarrow$ $t<\tau$:\\
This statement is proved in the same way as (ii).
\par
{\bf (iv)} So we get for each $\hat x\in\R^n$ numbers
$t_1,t_2\in[-\infty,+\infty]$ with $t_1\le t_2$ such that
$$\{\hat x\}\times(-\infty,t_1)\subset C,$$
$$\{\hat x\}\times([t_1,t_2]\cap\R)\subset\partial C$$
and
$$\{\hat x\}\times(t_2,+\infty)\subset\complement C.$$
$H^n$-almost eyerywhere in $\R^n$ we have
$t_1(\hat x)=t_2(\hat x)$: Show only that $H^n$-almost
everywhere the inequality $t_2(\hat x)-t_1(\hat x)<\frac{1}{k}$
for $k\ge1$ is valid.
Define $A_k:=\{\hat x\in\R^n:t_2(\hat x)-t_1(\hat x)\ge\frac{1}{k}\}$.
We get $0=H^{n+1}(\partial C)\ge\frac{1}{k}H^n(A_k)$ and
therefore $H^n(A_k)=0.$
Define $u(\hat x):=t_1(\hat x)$. We get
$C=\sub u$ in $L^1_{loc}(\rone{})$ and $u$ is measurable as
$C$ is measurable.
\qed

\def\ea{{\eta^\alpha}}
\def\eb{{\eta^\beta}}
\def\eat{{\tilde{\eta}^\alpha}}
\def\ebt{{\tilde{\eta}^\beta}}
\def\eae{{\eta^\alpha_\varepsilon}}
\def\ebe{{\eta^\beta_\varepsilon}}
\def\eate{{\tilde{\eta}^\alpha_\varepsilon}}
\def\ebte{{\tilde{\eta}^\beta_\varepsilon}}

\begin{lem}[C.~Gerhardt, according to a Seminar]%
\label{fastminimalbleibterhaltenlemma}
Let $E$ be an almost minimal set in
$\Omega\Subset S^n\times\R=\tN=(\tN,\gab)$.
Regard $\tN$ as a subset of $\rone{}$.
Then $E$ is also almost minimal in
$\Omega\subset\rone{}=(\rone{},\delta_{\alpha\beta})$.
The constant $\lambda$ in the definition of almost
minimal remains unchanged for $0<\lambda\le\frac{1}{2}$.
\end{lem}
\bew 
{\bf (i)} Almost minimal in a manifold is defined
in the same way as in $\rone{}$. Now the balls are
geodesic balls and the perimeter is defined
using the divergence in the manifold.
Since geodesic and Euclidean balls are comparable, i.~e.\
in any relatively compact subset $A$ and for any $\rho$
such that $0<\rho<R(A)$ there exists a constant $c=c(A)>0$
such that $B^\tN_{c\rho}(x)\subset B^\rone_\rho(x)
\subset B^\tN_{\frac{1}{c}\rho}(x)$ for any geodesic ball
in the respective manifold with center $x\in A$.
Therefore it suffices to prove the statement for
Euclidean balls.
From now on we assume that the chosen coordinates are
such that the standard metric of $\rone{}$ is represented
by the metric tensor
$(\delta_{\alpha\beta})_{0\le\alpha,\beta\le n}$.
\par
{\bf (ii)} Let $F$ be an arbitrary Caccioppoli set in $\tN$.
Let $\varepsilon>0$, $x_0\in A$ be arbitrary. Choose
$\eta^\alpha_\varepsilon \in C^1_c(B_\rho(x_0))$,
${0\le\alpha\le n}$, such that
\def\eae{{\eta^\alpha_\varepsilon}}
\def\ebe{{\eta^\beta_\varepsilon}}
$\gab\eae\ebe\le1$ and
$$\int\limits_{B_\rho(x_0)\subset\tN}|D\chi_F|\le\varepsilon+
  \int\limits_{B_\rho(x_0)\subset\tN}\chi_FD_\alpha\eae.$$
It follows that
\begin{align*}
\int\limits_{B_\rho(x_0)\subset\tN}|D\chi_F|
\le&\varepsilon+\int\limits_{B_\rho(x_0)\subset\tN}\chi_FD_\alpha\eae\\
=&\varepsilon+\int\limits_{B_\rho(x_0)}\chi_F
 \frac{\partial}{\partial x^\alpha}(\sqrt{\overline{g}(x)}\eae).
\end{align*}
For $x\in B_\rho(x_0)$, we get for sufficiently small $R=R(A)$
and $\rho$ such that $0<\rho<R$
\begin{align*}
\gab(x_0)\frac{\sqrt{\g(x)}}{\sqrt{\g(x_0)}}\eae(x)
  \frac{\sqrt{\g(x)}}{\sqrt{\g(x_0)}}\ebe(x)
\le&\sup\limits_{y\in A}\sup\limits_{z\in B_\rho(y)}
  \left(\frac{\sqrt{\g(z)}}{\sqrt{\g(x_0)}}\right)^2\cdot\\
&\cdot((\gab(x_0)-\gab(x))\eta^\alpha(x)\eta^\beta(x)\\
&\quad+\gab(x)\eta^\alpha(x)\eta^\beta(x))\\
\end{align*}
where the right-hand side can be estimated from above by
$$1+\rho c(A).$$
These estimates depend especially on the mean value theorem.
By $c(A)$ we denote a constant which depends only on $A$ and
may change its value from line to line. We have especially
$R=R(A)=c(A)$.
If we define $\eate=\frac{\sqrt{\g(x)}}{\sqrt{\g(x_0)}}\eae(x)$,
we deduce
\begin{align*}
\int\limits_{B_\rho(x_0)\subset\tN}|D\chi_F|
\le&\varepsilon+\int\limits_{B_\rho(x_0)}\chi_F
  \frac{\partial}{\partial x^\alpha}(\sqrt{\g(x_0)}\eate)
  \displaybreak[1]\\
\le&\varepsilon+(1+\rho c(A))\cdot\sup\left\{
  \int\limits_{B_\rho(x_0)}\chi_F
  \frac{\partial}{\partial x^\alpha}
  \left(\sqrt{\g(x_0)}\ea\right):\right.\\
&\eta^\alpha\in
  C^1_c(B_\rho(x_0)),{0\le\alpha\le n},\\
&\left.\rule{0em}{4.1ex}
  \gab(x_0)\ea(x)\eb(x)\le1\mbox{ for } x\in B_\rho(x_0)\right\}
  \displaybreak[1]\\
\equiv&\varepsilon+(1+\rho c(A))\cdot
  P_{(\rone{},\g(x_0))}(F,B_\rho(x_0)).
\end{align*}
As $\varepsilon>0$ was an arbitrary number, we get
$$\int\limits_{B_\rho(x_0)\subset\tN}|D\chi_F|
  \le(1+\rho c(A))\cdot P_{(\rone{},\g(x_0))}(F,B_\rho(x_0)).$$
\par
{\bf (iv)} Let $F$ be again an arbitrary Caccioppoli set in $\tN$.
Let $0<\rho<R$, $x_0\in A$ and
$\eta^\alpha\in C^1_c(B_\rho(x_0))$ such that
$\gab(x_0)\ea(x)\eb(x)\le1$.
For $\eat(x)=\frac{\sqrt{\g(x_0)}}{\sqrt{\g(x)}}\ea(x)$, we deduce
as above
$$\gab(x)\eat(x)\ebt(x)\le1+\rho c(A),$$
and we infer
\begin{align*}
\int\limits_{B_\rho(x_0)\subset\tN}|D\chi_F|
\ge&\frac{1}{1+\rho c(A)}\int\limits_{B_\rho(x_0)}
  \chi_F\frac{\partial}{\partial x^\alpha}(\sqrt{\g(x)}\eat(x))\\
\ge&\frac{1}{1+\rho c(A)}\int\limits_{B_\rho(x_0)}
  \chi_F\frac{\partial}{\partial x^\alpha}(\sqrt{\g(x_0)}\ea(x)).\\
\end{align*}
Now we take the supremum in this inequality over all
$\ea\in C^1_c(B_\rho(x_0))$
such that $\gab(x_0)\ea(x)\eb(x)\le1$.
This yields
$$\int\limits_{B_\rho(x_0)\subset\tN}|D\chi_F|
  \ge\frac{1}{1+\rho c(A)}P_{(\rone{},\g(x_0))}(F,B_\rho(x_0)).$$
\par
{\bf (v)} The conformal equivalence means
$\gab(x_0)=\vartheta(x_0)\delta_{\alpha\beta}(x_0)$
with $0<\vartheta(x_0)<\infty$. We have
$$P_{(\rone{},\g(x_0))}(F,B_\rho(x_0))=
 \varphi(\vartheta(x_0))\int\limits_{B_\rho(x_0)}|D\chi_F|,$$
where $\varphi$ is a positive and continuous function.
\par
{\bf (vi)} Now assume again that $E\vartriangle F\Subset B_\rho(x_0)$.
Using the result above, the inequalities of (iii) and (iv) can be
written in the following form:
$$\int\limits_{B_\rho(x_0)\subset\tN}|D\chi_F|\le
  (1+\rho c(A))
  \varphi(\vartheta(x_0))\int\limits_{B_\rho(x_0)}|D\chi_F|$$
and
$$\int\limits_{B_\rho(x_0)\subset\tN}|D\chi_E|\ge
  \frac{1}{1+\rho c(A)}\varphi(\vartheta(x_0))
  \int\limits_{B_\rho(x_0)}|D\chi_E|.$$
These inequalities are valid for any $x_0\in A$.
From (ii) we deduce the estimate
$$\int\limits_{B_\rho(x_0)\subset\tN}|D\chi_E|\le
  \int\limits_{B_\rho(x_0)\subset\tN}|D\chi_F|+
  K\rho^{n+2\lambda}.$$
Since $\varphi(\vartheta(x_0))$
is a positive and continuous function of
$x_0$, combining the inequalities above yields
$$\int\limits_{B_\rho(x_0)}|D\chi_E|\le
  (1+\rho\tilde c(A))\int\limits_{B_\rho(x_0)}|D\chi_F|
  +(1+\rho c(A))c(A)\rho^{n+2\lambda}.\qquad (\ast)$$
Choose $r$ with $0<r<\rho$ such that
$$\int\limits_{B_\rho(x_0)\setminus B_r(x_0)}|D\chi_E|
  \le\omega_{n+1}\rho^n.$$
We get the estimate
\begin{align*}
\int\limits_{B_\rho(x_0)}|D\chi_{E\cup B_r(x_0)}|
\le&\int\limits_{B_\rho(x_0)\setminus B_r(x_0)}|D\chi_E|+
  H^n(\partial B_r(x_0))\\
\le&\omega_{n+1}\rho^n+(n+1)\omega_{n+1}r^n\\
\le&(n+2)\omega_{n+1}\rho^n=c(A)\rho^n.\\
\end{align*}
Choosing $F=E\cup B_r(x_0)$, we obtain
$$\int\limits_{B_\rho(x_0)}|D\chi_E|\le
  c(A)(\rho^n+\rho^{n+2\lambda})\le c(A)\rho^n.$$
We multiply this inequality with $\rho\tilde{c}(A)$,
add it to ($\ast$) and get
$$(1+\rho \tilde{c}(A))\int\limits_{B_\rho(x_0)}|D\chi_E|\le
  (1+\rho \tilde{c}(A))\int\limits_{B_\rho(x_0)}|D\chi_F|
  +c(A)(\rho^{n+2\lambda}+\rho^{n+1}).$$
Now replace $\lambda$ by $\min\{\lambda,\frac{1}{2}\}$ and conclude
$$\int\limits_{B_\rho(x_0)}|D\chi_E|
  \le\int\limits_{B_\rho(x_0)}|D\chi_F|
  +c(A)\rho^{n+2\lambda}.$$
This inequality states that $E$ is
almost minimal in the set $\Omega\subset\rone{}$
equipped with the Euclidean metric.
\qed

Most of the proof of Theorem \ref{hyperflaechenbisdimensionacht}
in \cite{cg} is independent of the dimension,
so the following proof contains mainly the
necessary changes to get a proof which is valid up to
$n=7$:
\par
\bewvon{Theorem \ref{hyperflaechenbisdimensionacht}}
As in the cited paper we get
$u_k\in C^{2,\alpha}(S^n)$, $k\ge3$, such that
$$\left.H\right|_{{\rm graph\,}u_k}=f-\gamma e^{-\mu u_k}[u_k-u_{k-1}]$$
and $u_1\le u_k\le u_{k-1}$.
These functions converge pointwise, since for a fixed $x\in S^n$,
we have a monotone decreasing sequence $u_k(x)$ which is
bounded from below. Define $u$ to be the pointwise limit of $u_k$ and
further $\varphi_k=\log u_k$ and $\varphi=\log u$.
As in \cite{cg}, we deduce that
$E:={\rm sub\,}\varphi=\{(x,t):t<\varphi(x),\,x\in S^n\}$
is almost minimal in the metric product $S^n\times\R$.
Using Lemmata \ref{randinvarianzlemma} and
\ref{fastminimalbleibterhaltenlemma} we can regard $E$
an an almost minimal subset in $\rone{}$. Now Lemma
\ref{alssubgraphdarstellbarlemma} yields that the blow-up cone $C$
around a point $x\in\partial E$ is representable as a subgraph,
thus $C$ is a directed cone,
and the regularity theorem implies that $C$ is a half space.
Therefore we have $x\in\partial^\star E$ for any $x\in\partial E$.
From now on the proof of the theorem in \cite{cg} is independent
of the dimension and can be used to deduce the final result.
\qed

\section{Splitting Theorem}\label{splittingkap}

\begin{lem}\label{directionremainslemma}
Let $M\times\R\subset\R^n\times\R$ be a measurable set which
is directed with respect to $(\hat x,t)$.
If $\hat x\not=0$, then $M$ is directed with respect to $\hat x$.
\end{lem}
\bew
Since $M\times\R$ is directed with respect to $(\hat x,t)$, we
deduce in view of \cite[Theorem 1.3.6]{ziemer}
$$\int\limits_{\hat{B}_\rho(\hat z)\times[a,b]}\chi_{M\times\R}\ge
  \int\limits_{\hat{B}_\rho(\hat z+\sigma\hat x)
  \times[a+\sigma t,b+\sigma t]}\chi_{M\times\R}$$
for $\hat z\in\R^n$, $\rho>0$, $-\infty<a<b<+\infty$, $\sigma>0$
and $\hat B_\rho(\hat z)=\{\hat x\in\R^n:|\hat x-\hat z|<\rho\}$.
This inequality implies
$$\int\limits_{\hat{B}_\rho(\hat z)}\chi_{M}\ge
  \int\limits_{\hat{B}_\rho(\hat z+\sigma\hat x)}\chi_{M},$$
because we have
$$(b-a)\int\limits_{\hat{B}_\rho(\hat z)}\chi_{M}=
  \int\limits_{\hat{B}_\rho(\hat z)\times[a,b]}\chi_{M\times\R}$$
and an equality of the same kind is valid
for the other integrals above.
Since $\hat z\in\R^n$, $\rho>0$ and $\sigma>0$ were chosen
arbitrarily, the statement follows.
\qed

The following splitting theorem contains only the case of a
singular cone. For regular cones, however, there is no need
of such a theorem, because regular cones are half-spaces.

\begin{thm}[Splitting Theorem for directed minimal cones]
\hfill\\
Let $C\subset\rone$ be a singular minimal cone which has
$k$ linearly independent directions. Then there is a singular
minimal cone $C_0\subset\R^{n+1-k}$ such that $C=C_0\times\R^k$
after a suitable rotation and translation of $C$.
\end{thm}
\bew
We may assume that the vertex of the cone and the origin
coincide. In view of Lemma \ref{directionremainslemma}
we have to prove the Theorem only for $k=1$. Then the
statement for $k>1$
follows by induction. Assume that $C$ is directed
with respect to $e_{n+1}$.

Define $u:\R^n\to[-\infty,+\infty]$ as in \cite{mir}
to be a measurable function such that $C=\sub u$.
Then $u$ is positive homogeneous of degree 1, because $C$ is
a cone.
According to \cite[Theorem 15.5]{giusti} the set
$P:=\{\hat x\in\R^n:u(\hat x)=+\infty\}$ is a minimal set.
Since $u$ is positive homogeneous of degree 1 we deduce that
$P$ is a minimal cone with its vertex at the origin. Define
also $N:=\{\hat x\in\R^n:u(\hat x)=-\infty\}$. $N$ is also
a minimal cone with its vertex at the origin.

The case $N=P=\emptyset$ cannot occur:
$N=P=\emptyset$ implies $u\in L^\infty_{loc}(\R^n)$
(\cite[Proposition 16.7]{giusti}). According to
\cite[Theorem~1, p.~317]{gia1} we deduce $u\in BV(\R^n)$.
Finally \cite{cg1} implies $u\in C^{0,1}(\R^n)$ and
$u\in C^\omega(\R^n)$ follows. This would imply that
$C$ is regular, a contradiction.

According to Lemma \ref{negkomplementlemma}
we can assume that $P\not=\emptyset$, because in the case
$P\not=\emptyset$ we replace $C$ by $-\complement C$ and
proceed in the same way as we do now.

We also remark that $P$ is different from a half-space;
otherwise we get $P\times\R\subset C$ and in view of
\cite[Theorem 15.5]{giusti} $P\times\R=C$, a contradiction
to the fact that $C$ is a singular cone. In the same
way it can be shown that $N$ is different from a half-space.

$C\not=\rone$ implies
$P\not=\R^n$. According to the definition of
$P$ we deduce that $P\times\R\subset C$.
Then \cite[Example 16.2]{giusti} implies that $P\times\R$ is
minimal and therefore a singular minimal cone.
Now \cite[Theorem 2.4]{omc}, a maximum principle for minimal
cones, states, that two minimal cones are equal if they
have the same vertex and one of them contains the other one.
We apply this theorem and get
$C=P\times\R$. Define $C_0:=P$ and the statement follows.
\qed

\begin{cor}
Let $C\subset\rone$ be a singular minimal cone with its
vertex at the origin. Then we have, after a suitable rotation,
$C=C_0\times\R^k$, where $k\in\N$ is the number of 
linearly independent directions
and $C_0$ is a singular minimal
non-directed cone. In the exeptional case, $k=0$, we have of course
$C=C_0$.
\end{cor}

The splitting theorem contains the
Regularity Theorem \ref{regsatz}.

\begin{cor}
Let $C\subset\rone$ be a nontrivial minimal cone with
its vertex at the origin. Suppose $C$ has $k$ linearly
independent directions, $k\in\N$. Then $C$ is a half-space
provided that $n+1-k\le7$ or
$H^{k}(\partial C\setminus\partial^\star C)=0$.
\end{cor}
\bew
Assume the statement were false and let $C\subset\rone$ be
a counterexample. Then the splitting theorem yields (after
a suitable rotation of $C$) that $C=C_0\times\R^k$, where $C_0$
is a singular minimal cone in $\R^{n+1-k}$.

Since $n+1-k\le7$ contradicts the non-existence of singular
minimal cones up to $\R^7$ (\cite{giusti}),
this case does not occur.

In the second case we deduce that
$H^0(\partial C_0\setminus\partial^\star C_0)\ge1$, because
$C_0$ is a singular minimal cone. This implies
\begin{align*}
H^k(\partial(C_0\times\R^k)\setminus\partial^\star(C_0\times\R^k))
=&H^k((\partial C_0\setminus\partial^\star C_0)\times\R^k)\\
\ge&H^0(\partial C_0\setminus\partial^\star C_0)\ge1\\
>&0=H^k(\partial C\setminus\partial^\star C),
\end{align*}
contradicting $C=C_0\times\R^k$.
Thus the statement is proved.
\qed

\begin{cor}\label{onepointcor}
Let $C\subset\rone$ be an open minimal cone with its vertex at
the origin. If $C$ is nontrivial and directed with respect to
$e_{n+1}$ and there exists a $t>0$ such that $(0,t)\in\complement C$
or $(0,-t)\in C$, then $\partial C$ is a hyperplane.
\end{cor}
\bew
In view of Lemma \ref{negkomplementlemma} we can assume that
$(0,-t)\in C$. So there exists $\rho>0$ such that
$B_\rho((0,-t))\subset C$.
Thus $u(\hat x)>-t$ holds for $|\hat x|<\rho$,
where $u$ is defined as in Remark \ref{subgleichger}
such that $\sub u=C$.
In view of the homogeneity of $u$, this implies
$u(\hat x)>-\infty$ for $\hat x\in\R^n$.
Therefore we deduce that
$N=\emptyset$, where the notation from the proof of the
splitting theorem has been used and will be used in the rest
of this proof. Since $P=\emptyset$ implies the statement as
shown in the proof of the splitting theorem, we may assume assume
$\emptyset\not=P\not=\R^n$. Following again the proof of the
splitting theorem, we deduce that $C=P\times\R$. We choose now
$\hat x\in\complement P$. Since $\complement P$ is a cone,
we have $\tau\hat x\in\complement P$ for $\tau>0$ and
$(\tau\hat x,-t)\in\complement C$, but $(\tau\hat x,-t)$
converges to $(0,-t)\in C$ for $\tau\to0$. This is a
contradiction, because $C$ is an open set.
\qed

\begin{rem}
The assumptions of Corollary \ref{onepointcor} imply
also $\langle\nu,e_{n+1}\rangle<0$, where $\nu$ is the inner
unit normal of $C$.
\end{rem}

\bibliographystyle{amsplain}

\begin{thebibliography}{99}
\bibitem{dg} E.~De~Giorgi: {\it Una estensione del teorema
  di Bernstein,\/} Ann.\ Scuola Norm.\ Sup.\ Pisa (3) {\bf 19}
  (1963), 79-85.
\bibitem{cg} C.~Gerhardt: {\it Closed Hypersurfaces of prescribed
  Mean Curvature in locally conformally flat Riemannian Manifolds,\/}
  J.~Diff.~Geom.~{\bf 48} (1998), 587-613.
\bibitem{cg1} C.~Gerhardt: {\it On the Regularity of Solutions to
  Variational Problems in $BV(\Omega)$,\/} Math.~Z.~{\bf 149} (1976),
  281-286.
\bibitem{gia1} M.~Giaquinta, G.~Modica, J.~Sou\v{c}ek:
  {\it Cartesian Currents in the Calculus of Variations I,\/}
  Ergebnisse der Mathematik und ihrer Grenzgebiete. 3.~Folge,
  Band~37, Springer, Berlin Heidelberg New~York, 1998, 736~pages.
\bibitem{giusti}
  E.~Giusti: {\it Minimal Surfaces and Functions of Bounded Variation,\/}
  Monographs in Math.\ Vol.~80, Birkh\"auser, Boston Basel Stuttgart,
  1984, 240~pages.
\bibitem{omc}
  E.~Gonzalez, U.~Massari, M.~Miranda:
  {\it On minimal Cones,\/} Appl.\ Anal.~{\bf 65} (1997), 135-143.
\bibitem{mm} U.~Massari, M.~Miranda:
  {\it Minimal Surfaces of Codimension One,\/} Mathmatics Studies Vol.\
  91, North-Holland, Amsterdam New~York Oxford, 1984, 242~pages.
\bibitem{mir} M.~Miranda:
  {\it Superficie minime illimitate,\/}
  Ann.\ Sc.\ Norm.\ Sup.\ Pisa (4) {\bf 4} (1977), 313-322.
\bibitem{tam2} I.~Tamanini:
  {\it Boundaries of Caccioppoli sets with H\"older-continuous
  normal vector,\/} J. Reine Angew.\ Math.\ {\bf 334} (1982), 27-39.
\bibitem{tamanini} I.~Tamanini:
  {\it Regularity results for almost minimal oriented hypersurfaces
  in $\R^N$,\/} Quaderni Dipartimento Math.\ Univ.\ Leece {\bf 1},
  (1984), 92~pages.
\bibitem{ziemer} W.~P.~Ziemer: {\it Weakly differentiable Functions,\/}
  Graduate Texts in Mathmatics Vol. 120, Springer,
  Berlin Heidelberg New~York, 1989, 308~pages.
\end{thebibliography}

\end{document}